\documentclass{amsart}
\usepackage{amssymb}
\usepackage{mathrsfs}
\usepackage{amscd}

\let\mathcal \undefined
\def\mathcal{\mathscr}

\usepackage[colorlinks,linkcolor={blue},
citecolor={blue},urlcolor={red},]{hyperref}

\theoremstyle{plain}
\newtheorem{theorem}{Theorem}
\theoremstyle{remark}

\theoremstyle{plain}

\theoremstyle{plain}

\def\R{{\mathbb R}\,}

\def\E{{\mathbb E}\,}

\renewcommand{\ell}{\g}
\newcommand{\one}{{{\bf 1}}}

\newcommand{\calL}{{\mathcal L}}

\newcommand{\g}{\gamma}

\newcommand{\n}{\Vert}

\newcommand{\e}{\varepsilon}

\newcommand{\Lipo}{{\rm Lip}_0(X,Y)}

\renewcommand{\l}{l_n^2}

\begin{document}
\title
{On the action of Lipschitz functions on vector-valued random sums}

\author{Jan van Neerven}
\author{Mark Veraar}

\address{Delft Institute of Applied Mathematics\\
Technical University of Delft \\ P.O. Box 5031\\ 2600 GA Delft\\The
Netherlands}

\email{J.vanNeerven@math.tudelft.nl, M.C.Veraar@math.tudelft.nl}

\thanks{The authors are supported by the `VIDI subsidie' 639.032.201
in the `Vernieuwingsimpuls' programme of the Netherlands
Organization for Scientific Research (NWO) and
by the Research Training Network HPRN-CT-2002-00281.}

\keywords{Lipschitz functions, type 2, cotype 2,
isomorphic characterization of Hilbert spaces, Dvoretzky's theorem}

\subjclass[2000]{Primary: 46C15, Secondary: 46B09, 47B10}

\date\today

\begin{abstract}
Let $X$ be a Banach space and let $(\xi_j)_{j\ge 1}$ be an i.i.d. sequence of
symmetric
random variables with finite moments of all orders.
We prove that the following assertions are equivalent:
\begin{enumerate}
\item There exists a constant $K$
such that
$$
\Bigl(\E\Big\|\sum_{j=1}^n \xi_j f(x_j)\Big\|^2\Bigr)^{\frac12}
\leq K \n f\n_{\rm Lip}
\Bigl(\E\Big\|\sum_{j=1}^n \xi_j x_j\Big\|^2\Bigr)^{\frac12}
$$
for all Lipschitz functions $f:X\to
X$ satisfying $f(0)=0$ and all finite sequences $x_1,\dots,x_n$ in $X$.
\item $X$ is isomorphic to a Hilbert space.
\end{enumerate}
\end{abstract}

\maketitle

For Banach spaces $X$ and $Y$ let $\Lipo$ denote the Banach
space of all Lipschitz continuous functions $f:X\to Y$ satisfying $f(0)=0$
with norm $\n f\n_{\rm Lip}:= L_f,$
the Lipschitz constant of $f$.
Our main result relates the action of functions $f\in\Lipo$
on random sums in $X$ with the cotype and type of $X$ and $Y$, respectively.
Since the best constants are obtained for Gaussian variables, we state the
result for this case first.

\begin{theorem}\label{thm:main}
Let $X$ and $Y$ be Banach spaces with $\dim X=\infty$ and $\dim Y\ge 1$,
and let $(\g_j)_{j\ge 1}$ be a sequence of independent standard Gaussian random
variables.
The following assertions are equivalent:
\begin{enumerate}
\item[(i)] For all finite sequences
$x_1, \ldots, x_n\in X$, all scalars $a_1,\dots,a_n >0$,
and all $f_1,\dots,f_n \in \Lipo$ we have
\[\Bigl(\E\Big\|\sum_{j=1}^n \g_j a_j^{-1} f_j(a_j x_j)\Big\|^2\Bigr)^{\frac12}
\leq K \big(\max_{1\le j\le n}\n f_j\n_{\rm Lip}\big) \Bigl(\E\Big\|\sum_{j=1}^n \g_j x_j\Big\|^2\Bigr)^{\frac12},\]
where $K$ is a constant depending on $X$ and $Y$ only.

\item[(ii)] For all finite sequences
$x_1, \ldots, x_n\in X$ there exist scalars $a_1,\dots,a_n>0$ such that for all
$f\in \Lipo$ we have
\[\Bigl(\E\Big\|\sum_{j=1}^n \g_j a_j^{-1} f(a_j x_j)\Big\|^2\Bigr)^{\frac12}\leq K \n f\n_{\rm Lip} \Bigl(\E\Big\|\sum_{j=1}^n \g_j x_j\Big\|^2\Bigr)^{\frac12},\]
where $K$ is a constant depending on $X$ and $Y$ only.

\item[(iii)] $X$ has cotype $2$ and $Y$ has type $2$.
\end{enumerate}
If ${\rm (i)}$ or ${\rm (ii)}$ holds with constant $K$,
then the Gaussian cotype $2$ constant of $X$ and
the Gaussian type $2$ constant of $Y$ satisfy
$C_2^\g(X)\le K$ and $T_2^\g(Y)\le \sqrt{2}K.$
\end{theorem}

\begin{proof}
The implication (i)$\Rightarrow$(ii) is trivial.

(ii)$\Rightarrow$(iii):
First
we prove that $X$ has cotype $2$ with $C_2^\g(X)\le K$.
Fix a norm one vector $y_0\in Y$ and
define $f\in\Lipo$ by $f(x) := \n x\n y_0.$ Since $\n f\n_{\rm Lip} = 1$ it
follows that for $x_1, \ldots, x_n\in X$ we have, with the
$a_1,\dots,a_n>0$ as in (ii),
$$
\sum_{j=1}^n \n x_j\n^2  = \E\Big\|\sum_{j=1}^n \g_j f(x_j)\Bigr\|^2  =
\E\Big\|\sum_{j=1}^n \g_j a_j^{-1} f(a_j x_j)\Bigr\|^2 \leq K^2 \E
\Big\n \sum_{j=1}^n \g_j \, x_j\Big\n^2.
$$

Next we prove that $Y$ has type $2$ with
$T_2^\g(Y)\le K\sqrt{2}$.
By an observation in \cite{Ja} we have
\begin{equation}\label{pisier}
T_2^\g(Y) = \sup\Big\{\E \Big(\Big\n\sum_{j=1}^n \g_j \,y_j\Big\n^2
\Bigr)^\frac12:\ n\ge 1, \ \n
y_1\n = \dots = \n y_n\n = n^{-\frac12}\Big\}.
\end{equation}
Fix an integer $n\ge 1$ and  vectors $y_1,\dots,y_n\in Y$ of norm $1$. Let
$(e_j)_{j=1}^n$ be the standard unit basis of $\l$ and let $\e>0$ be arbitrary
and fixed. Since $\dim X = \infty$, by Dvoretzky's theorem \cite{Dv} we can
find an isomorphism $T$ from $\l$ onto an $n$-dimensional subspace $X_0$ of
$X$ such that
 $\|T\|\le 1+\e$ and $\n T^{-1}\n = 1$. Let
 \begin{equation}\label{defxj}
 x_j := T e_j, \ \  j=1,\dots,n.
 \end{equation}
 Clearly, $1\le\n x_j\n\le 1+\e$ and for all $1\le j\neq k\le n$
and $a, b\in \R$ we have
\begin{equation}\label{sqrt}
\begin{aligned}
\n a x_j- b x_k\n & \ge \n T^{-1}\n^{-1}\n a e_j-b e_k\n
 = \sqrt{a^2+b^2}.\end{aligned}
\end{equation}

Define $\varphi_j:X\to \R$ by
$$\varphi_j(x):= \max\big\{0, 1-\sqrt{2}\n x-x_j\n\big\}.$$
Then $\varphi_j$ is Lipschitz continuous with Lipschitz constant
$\n{\varphi_j}\n_{\rm Lip}\le \sqrt{2}$, we have $\varphi_j(x_j)=1$,
and $\varphi_j\equiv 0$
outside the open `sector'
\[V_j := \Big\{x\in X: \exists t>0 \ \text{such that} \ \|t x - x_j\|
< \tfrac12\sqrt{2}\Big\}.\] Note that $0\notin V_j$. We claim that the sectors
$V_j$ are disjoint. Indeed, given $x\in V_j$ we choose $t>0$ such that $\|t x
- x_j\|<\frac12\sqrt{2}$. Then for $j\neq k$ and all $s>0$,
\[\|s x- x_k\| \geq \|t^{-1} s x_j - x_k\| - \|t^{-1} s x_j - s x\|
\stackrel{{\rm (*)}}{>} \sqrt{t^{-2} s^2 + 1} - \tfrac12  t^{-1}
s\sqrt{2}\stackrel{{\rm (**)}}{\geq} \tfrac12 \sqrt{2}.\] In $(*)$ we used
\eqref{sqrt} and the choice of $t$, while $(**)$ follows from the inequality
$\sqrt{c^2 + 1} - \frac12 c\sqrt{2}\geq \frac12 \sqrt{2}$.

Define $\psi_j: X\to \R$ by
\[\psi_j(x) := a_j \varphi_j\big(a_j^{-1} x\big),\]
where the $a_1,\dots,a_n>0$ are chosen as in (ii).
Then $\psi_j$ is Lipschitz continuous with Lipschitz constant
$\n{\psi_j}\n_{\rm Lip}\leq \sqrt{2}$, we have $\psi_j(a_j x_j) =
a_j$, and $\psi_j\equiv 0$ outside $V_j$. Define $f:X\to Y$ by
$$ f(x) := \sum_{j=1}^n \psi_j(x)y_j.$$
It is clear that $f(0)=0$ and $f(a_j x_j) = a_jy_j$.
We claim that $f\in \Lipo$ with $\n f\n_{\rm Lip} \leq \sqrt{2}$. If $x,x'\in V_j$ for some
$j$, then by the disjointness of $V_j$ with the other $V_k$'s and the fact
that $\n y_j\n=1$ we obtain
$$ \n f(x)-f(x')\n = \n y_j\n \, |\psi_j(x)-\psi_j(x')| \le \sqrt{2}\,\n x-x'\n.$$
If $x\in V_j$ and $x'\in V_k$ for $j\not=k$, we choose convex combinations
$\xi$ and $\xi'$ of $x$ and $x'$, say $\xi = (1-s)x+sx'$ and $\xi'=(1-t)x+tx'$
with $0\leq s\leq t\leq 1$, such that $\xi\in\partial V_j$ and
$\xi'\in\partial V_k$. Clearly, $f(\xi)=f(\xi')=0$. It follows from the
previous case that
\[\begin{aligned}\|f(x)-f(x')\|&\leq \|f(x)-f(\xi)\| +
\|f(\xi')-f(x')\|
\\ &\leq \sqrt{2} \|x-\xi\|+\sqrt{2}\|\xi'-x'\|
\\ & = \sqrt{2}(s+(1-t))\|x-x'\|\leq \sqrt{2}\|x-x'\|.
\end{aligned}
\]
The case where $x\in V_j$ and $x' \notin \bigcup_{k}V_k$ is handled similarly.
Finally if $x,x' \notin \bigcup_{k}V_k$, then $f(x) = f(x')=0$. This concludes
the proof of the claim.

Recalling that $f(0)=0$, $\n f\n_{\rm Lip}\le \sqrt{2}$, $\n T \n\le 1+\e$, $\n y_j\n=1$,
we obtain
$$
\begin{aligned}
 \E \Big\n \sum_{j=1}^n \g_j \, y_j\Big\n^2
& = \E \Big\n \sum_{j=1}^n \g_j \, a_j^{-1} f(a_j x_j) \Big\n^2 \leq
2K^2 \E \Big\n \sum_{j=1}^n \g_j \, x_j \Big\n^2
\\ &   \le   2 K^2(1+\e)^2 \E \Big\n \sum_{j=1}^n \g_j \, e_j\Big\n^2
 =  2 K^2(1+\e)^2\sum_{j=1}^n \n y_j\n^2.
 \end{aligned}
 $$
By \eqref{pisier} this proves that $Y$ has type $2$ with $T_2^\g(Y)\le K
\sqrt{2}(1+\e)$. Since $\e>0$ was arbitrary, the proof is complete.

(iii)$\Rightarrow$(i): Assume that $X$ has cotype $2$ and $Y$ has type $2$.
For all $x_1, \ldots, x_n\in X$,  $a_1,\dots,a_n>0$, and $f_1,\dots,f_n\in\Lipo$ we have
\[
\begin{aligned}
\E\Big\|\sum_{j=1}^n \g_j a_j^{-1} f_j(a_j x_j)\Big\|^2
& \leq T_2^\g(Y)^2\big(\max_{1\le j\le n} \n f_j\n_{\rm Lip}\big)^2 \sum_{j=1}^n a_j^{-2} \|a_j x_j\|^2
\\ & \leq T_2^\g(Y)^2\big(\max_{1\le j\le n} \n f_j\n_{\rm Lip}\big)^2 C_2^\g(X)^2 \E\Big\|\sum_{j=1}^n \g_j x_j\Big\|^2.
\end{aligned}
\]
\end{proof}

By a celebrated theorem of Kwapie\'n \cite{Kwa}, a Banach space $X$ has type
$2$ and cotype $2$ if and only if $X$ is isomorphic to a Hilbert space. Thus
if we take $X=Y$ in the theorem, then assertion (iii)  may be replaced by:
\smallskip

\begin{enumerate}
\item[(iii)$'$] $X$ is isomorphic to a Hilbert space.
\end{enumerate}
\smallskip

In Theorem \ref{thm:main} we may replace the Gaussian sequence $(\g_j)_{j\geq
1}$ by a Rademacher sequence $(r_j)_{j\geq 1}$, in which case
we obtain the estimates $$C_2^{r}(X)\le K\quad\text{ and }\quad
T_2^r(Y)\le \frac{2}{\sqrt{\pi}}K.$$
Here $C_2^{r}(X)$ and $T_2^{r}(Y)$ denote
the Rademacher cotype $2$ constant of
$X$ and the Rademacher type $2$ constant of
$Y$, respectively. For the second estimate we recall from
\cite[Lemma 4.5]{LeTa} that
$T_2^r(X)\leq \frac1{m_1^\g} T_2^\g(X)$,
where $m_1^\g:= \E|\g_j| = \sqrt{2/\pi}$ and
that by an observation in \cite{Ja} we have
\begin{equation}\label{pisier2}
T_2^\g(Y) = \sup\Big\{\E \Big(\Big\n\sum_{j=1}^n r_j \,y_j\Big\n^2\Bigr)^\frac12 :\ n\ge 1, \ \n y_1\n = \dots = \n y_n\n =
n^{-\frac12}\Big\}.\end{equation}
The proof of (ii)$\Rightarrow$(iii) may now be repeated verbatim.

For Banach spaces $X$ with the finite Lipschitz extension property it is possible to give a considerable simpler proof of Theorem 1.

Next let $(\xi_{j})_{j\geq 1}$ be an arbitrary sequence of i.i.d.
symmetric random variables with $\E|\xi_j|^2=1$.
We denote by $T_2^\xi(X)$ and $C_2^{\xi}(X)$ the $\xi$-type $2$ and
$\xi$-cotype $2$ constant of a Banach space, respectively.
By a standard randomization argument, every Banach space $X$ with (co)type $2$
has $\xi$-(co)type $2$ with constants $T_2^\xi(X)\leq T_2^r(X)$ and
$C_2^\xi(X)\leq C_2^r(X)$.
Conversely, if $X$ has $\xi$-type $2$, then again by \cite[Lemma
4.5]{LeTa},
\[\Big(\E\Big\|\sum_{j=1}^n r_j x_j\Bigr\|^2\Bigr)^\frac12
\leq \frac1{m_1^\xi}\Big(\E\Big\|\sum_{j=1}^n \xi_jx_j\Bigr\|^2\Bigr)^\frac12
\leq \frac1{m_1^\xi} T_2^\xi(X) \Big(\sum_{j=1}^n \n x_j\n^2\Bigr)^\frac12,\]
where $m_1^\xi:= \E|\xi_j|$. It follows that $X$ has type $2$ with
$T_2^r(X)\leq \frac1{m_1^\xi} T_2^\xi(X)$.
If $X$ has $\xi$-cotype $2$ and all moments of $\xi_j$ are finite,
then $X$ has finite cotype (we are grateful to
Tuomas Hyt\"onen for pointing
this out to us). In fact, by means of elementary estimates it can be shown that $c_0$
does not have finite $\xi$-cotype.
The Rademacher cotype $2$ of $X$ then follows from the Maurey-Pisier theorem; cf.
\cite[Section 9.2]{LeTa}.

At the expense of slightly worse estimate for the type $2$ constant
it is possible to
generalize Theorem \ref{thm:main} to sequences of random variables
$(\xi_j)_{j\ge 1}$ as above.
This is achieved by a slightly modified argument
which does not require normalizations as
in \eqref{pisier} and \eqref{pisier2} and which has the additional virtue
that for each $n$ the scalars $a_{1},\dots, a_{n}$ are allowed to
depend not only on the vectors $x_1,\ldots, x_n$ but also on the function $f$.

\begin{theorem}\label{thm:main2}
Let $X$ and $Y$ be Banach spaces with $\dim X=\infty$ and $\dim Y\ge 1$,
and let $\xi=(\xi_j)_{j\ge 1}$ be a sequence of i.i.d. random variables with
$\E|\xi_j|^2 = 1$. The following assertions
are equivalent:
\begin{enumerate}
\item[(i)] For all $f_1,\dots,f_n\in\Lipo$, all finite sequences
$x_1, \ldots, x_n\in X$, and all scalars $a_{1},\dots, a_{n}>0$ we have
\[\Bigl(\E\Big\|\sum_{j=1}^n \xi_j a_j^{-1} f_j(a_j x_j)\Big\|^2\Bigr)^{\frac12}\leq K \big(\max_{1\le j\le n} \n f_j\n_{\rm Lip}\big) \Bigl(\E\Big\|\sum_{j=1}^n \xi_j x_j\Big\|^2\Bigr)^{\frac12},\]
where $K$ is a constant depending on $X$ and $Y$ only.

\item[(ii)] For all  $f\in \Lipo$ and all
finite sequences $x_1, \ldots, x_n\in X$
there exist scalars $a_{1},\dots, a_{n}>0$ such that
\[\Bigl(\E\Big\|\sum_{j=1}^n \xi_j a_j^{-1} f(a_j x_j)\Big\|^2\Bigr)^{\frac12}\leq K \n f\n_{\rm Lip} \Bigl(\E\Big\|\sum_{j=1}^n \xi_j x_j\Big\|^2\Bigr)^{\frac12},\]
where $K$ is a constant depending on $X$ and $Y$ only.

\item[(iii)] $X$ has $\xi$-cotype $2$ and $Y$ has $\xi$-type $2$.
\end{enumerate}
If ${\rm (ii)}$ holds, then $C_2^{\xi}(X)\le K$ and $T_2^{\xi}(Y)\le
(1+2\sqrt{2})K$.
If the $\xi_j$ have finite moments of all orders, then {\rm (iii)} is equivalent to
\begin{enumerate}
\item[(iv)] $X$ has cotype $2$ and $Y$ has type $2$.
\end{enumerate}
\end{theorem}

\begin{proof}
Only the proof that $Y$ has $\xi$-type $2$ in the implication
(ii)$\Rightarrow$(iii) needs to be adapted.
Fix arbitrary nonzero vectors $y_1,\dots,y_n\in Y$.
Following the arguments in the proof of (ii)$\Rightarrow$(iii) in Theorem \ref{thm:main},
we replace \eqref{defxj} by
$$
x_j := \|y_j\| T e_j, \ \ j=1,\dots,n,$$
and define $\varphi_j:X\to \R$ by $\varphi(0)=0$ and
$$\varphi_j(x):= \max\Big\{0, 1-\sqrt{2} (1+\e) d_j(x)\Big\}\|x\| ,$$
where $d_j:X\setminus\{0\}\to \R$ is the function
\[d_j(x) :=\Big\n \frac{x}{\n x\n}- \frac{x_j}{\|x_j\|}\Big\n.\]
Then $\varphi_j$ is Lipschitz continuous with
$\|\varphi_j\|_{\rm Lip}\le L_{\e}:=2\sqrt{2}(1+\e)+1$,
we have $\varphi_j(a x_j)=a\|x_j\|$ for all  $a>0$, and
$\varphi_j\equiv 0$ outside the sector
\[V_j := \Big\{x\in X\setminus \{0\}:  \ d_j(tx)
< \tfrac12\sqrt{2}(1+\e)^{-1}\Big\}.\]
As before, $V_j$ and $V_k$ are disjoint for $j\neq k$.
Indeed if $x\in V_j$, then for $j\not=k$ we have
\[\begin{aligned}
\Big\|\frac{x}{\|x\|}- \frac{x_k}{\|x_k\|} \Big\| &\geq
\Big\|\frac{x_j}{\|x_j\|}- \frac{x_k}{\|x_k\|} \Big\| -
\Big\|\frac{x_j}{\|x_j\|}- \frac{x}{\|x\|} \Big\| \\ & > \sqrt{\|T
e_j\|^{-2}+\|T e_k\|^{-2}} - \tfrac12\sqrt{2}(1+\e)^{-1}
\\ & \geq \sqrt{2} (1+\e)^{-1} - \tfrac12\sqrt{2}(1+\e)^{-1} =
\tfrac12\sqrt{2}(1+\e)^{-1},
\end{aligned}\]
which shows that  $x\notin S_k$.
Define $f:X\to Y$ by
$$ f(x) = \sum_{j=1}^n \varphi_j(x)\frac{y_j}{\|x_j\|}.$$
Then $f(0)=0$, $f(a x_j) = a y_j=a f(x_j)$ for $a>0$, and  $f$ is Lipschitz continuous with
$\|f\|_{\rm Lip}\le L_{\e}$.
With the $a_1,\dots,a_n>0$ as in (ii), estimating as before we obtain
$$
\E \Big\n \sum_{j=1}^n \xi_j \, y_j\Big\n^2
= \E \Big\n \sum_{j=1}^n \xi_j \, a_{j}^{-1} f(a_{j}x_j)\Big\n^2
\le  \|f\|_{\rm Lip}^2 K^2(1+\e)^2\sum_{j=1}^n \n y_j\n^2.
$$
This proves that $Y$ has $\xi$-type $2$ with
\[T_2^{\xi}(Y)\le K \|f\|_{\rm Lip} (1+\e) \leq K \big(1+2\sqrt{2} (1+\e)\big) (1+\e).\]
Since $\e>0$ was arbitrary, the proof is complete.
\end{proof}

If the $\xi_j$ have finite moments of all orders, for $X=Y$ we obtain
an isomorphic characterization of Hilbert spaces as before.
\medskip

Theorems \ref{thm:main} and \ref{thm:main2} bear
a striking resemblance to \cite[Proposition 1.13]{ArBu}
which states that $X$ has type $2$ and $Y$ has
cotype $2$ if and only if every uniformly bounded family $\mathcal{T}$ in
$\calL(X,Y)$ is $R$-bounded. Recall that $\mathcal{T}$ is called
{\em $R$-bounded} if there exists a constant $K$ such that
for all choices $x_1,\dots,x_n\in X$ we have
$$ \Big(\E\Big\n \sum_{j=1}^n r_j\,T_j x_j\Big\n^2\Big)^\frac12
\le K \Big(\E\Big\n \sum_{j=1}^n r_j\,x_j\Big\n^2\Big)^\frac12.
$$
This result is elementary (it suffices to consider suitably chosen
families of rank one operators) and the role of
the Rademacher variables can be replaced by any i.i.d. sequence of mean zero
random variables with finite second moment.
The precise relationship between
\cite[Proposition 1.13]{ArBu} and our results remains unclear, since we see no obvious way to relate finitely many linear operators
in $\calL(X,Y)$ to a single nonlinear function in $\Lipo$.
In this connection it is worthwhile to point out that it appears to be an unsolved open problem
whether for every pair of Banach spaces $X$ and $Y$ there exists a constant $c(X,Y)$ such that, given any
distinct elements $x_1,\dots,x_n\in X$ and elements
$y_1,\dots,y_n\in Y$,  there exists a Lipschitz function
$f:X\to Y$ satisfying $f(x_j) = y_j$ for all $j=1,\dots,n$
and
\begin{equation}\label{interpol} \n f\n_{\rm Lip}
\le c(X,Y) \max_{\substack{1\le j,k\le n\\ j\not=k}}
\frac{\n y_j-y_k\n}{\n x_j-x_k\n}.
\end{equation}
The important point here is that $c(X,Y)$ should be independent of $n$.
Indeed, it was shown in \cite{JLS}
that for fixed $n$, \eqref{interpol} can be achieved with a constant
$c(n,X,Y)$ of order $\log n$.

\medskip
As an application of Theorem \ref{thm:main} we will prove next that
Lip$_0(X)$ acts in the operator ideal $\g(l^2,X)$ of
$\g$-radonifying operators from $l^2$ to $X$ if and only if $X$ is isomorphic to
a Hilbert space.

Let $H$ be a Hilbert space. We denote by $\ell(H,X)$
the completion of the vector space of all finite rank operators $u:H\to X$
with respect to the norm
\begin{equation}\label{gnorm}
\|u\|_{\ell(H,X)} := \sup \Big(\E \Big\n \sum_{j} \g_j \, u
h_j\Big\n^2\Big)^\frac12.
\end{equation}
The supremum is taken over all finite orthonormal systems $(h_j)$ in $H$. As
is well known, $\ell(H,X)$ is an operator ideal in the sense that for all
bounded linear operators $v:\tilde H\to H$ and $w:X\to \tilde X$ we have
$wuv\in\ell(\tilde H, \tilde X)$ and $$\n wuv\n_{\ell(\tilde H, \tilde X)}\le
\n w\n\, \n u\n_{\ell(H,X)}\,\n v\n.$$
For more information we refer to
\cite[Chapter 12]{DJT}.

We will be interested in the particular case where $H$ equals $L^2:=
L^2(S,\Sigma,\mu)$ for some $\sigma$-finite measure space $(S,\Sigma,\mu)$ and
$u_\phi: L^2\to X$ is an integral operator of the form
\[u_\phi h = \int_S h(s)\phi(s) \, d\mu(s), \ \ h\in L^2,\]
for suitable functions $\phi:S\to X$. Operators in $\ell(L^2,X)$ arising in
this way have been investigated recently in \cite{KW}. If $\phi$ is a simple
function, i.e., a function of the form $\sum_{j=1}^n \one_{S_j} \otimes x_j$
with vectors $x_j$ taken from $X$ and disjoint sets $S_j\in\Sigma$ satisfying
$0<\mu(S_j)<\infty$, it is easily checked that $u_\phi\in
\ell(L^2, X)$ and  by considering the orthonormal functions $h_j :=
\mu(S_j)^{-\frac12}\one_{S_j}$, the $\ell$-norm of $\phi$ is computed as
\begin{equation}\label{lnorm}
\|u_\phi\|_{\ell(L^2, X)}^2 = \E \Big\n \sum_{j=1}^n \g_j u_\phi h_j  \Big\n^2
= \E \Big\n \sum_{j=1}^n \g_j \, \mu(S_j)^{\frac12} x_j \Big\n^2.
\end{equation}
The subspace of all $u\in \ell(L^2,X)$ of the form $u=u_\phi$ for some simple
function $\phi:S\to X$ will be denoted by $\ell_{\text{simple}}(L^2,X)$. An easy
approximation argument shows that this is a dense subspace of $\ell(L^2,X)$.

If $X$ has type $2$, the mapping $\phi\mapsto u_\phi$ defined for simple
functions $\phi$ as above,
extends to a continuous embedding
from $L^2(X):=L^2(S,\Sigma,\mu;X)$ into $\ell(L^2, X)$.
Indeed, for a simple function $\phi = \sum_{j=1}^n \one_{S_j}\otimes x_j$ we have, using \eqref{lnorm},
\begin{equation}\label{type2}
\begin{aligned}
\|u_{\phi}\|_{\ell(L^2, X)}^2
& = \E \Big\n \sum_{j=1}^n \g_j \,
\mu(S_j)^{\frac12} x_j \Big\n^2
\\ & \leq T_2^\g(X)^2 \sum_{j=1}^n \mu(S_j) \n
x_j\n^2 = T_2^\g(X)^2 \|\phi\|_{L^2(X)}^2,
\end{aligned}
\end{equation}
and the claim follows by a density argument. Similarly, if $X$ has cotype $2$,
then $u_\phi\mapsto \phi$ extends to a continuous embedding from
$\ell(L^2, X)$ into $L^2(X)$.

If $\phi =\sum_{j=1}^n \one_{S_j}\otimes x_j$ is a simple $X$-valued function,
then for each $f\in\Lipo$,
\[f(\phi) = \sum_{j=1}^n \one_{S_j}\otimes f(x_j)\]
is a simple $Y$-valued function. In this way we obtain a mapping
$\tilde{f}:\ell_{\rm simple}(L^2, X) \to \ell_{\rm simple}(L^2, Y)$ by putting
\[\tilde{f}(u_{\phi}) := u_{f(\phi)}.\]
We are interested in conditions ensuring that $\tilde f$ extends to a
Lipschitz continuous mapping from $\ell(L^2, X)$ to $\ell(L^2, Y)$.
From $f(0)=0$ we see that a necessary condition
is that there should exist a constant $K$ such that
\[\|u_{f(\phi)}\|_{\ell(L^2, Y)}\leq K \n f\n_{\rm Lip} \|u_{\phi}\|_{\ell(L^2, X)}\]
for all simple functions $\phi:S\to X$. The next result gives a converse and
relates both conditions to the geometry of the spaces $X$ and $Y$.

\begin{theorem}\label{thm:main3}
Let $X$ and $Y$ be Banach spaces, let $L^2:=L^2(S,\Sigma,\mu)$ as before,
and assume that $\dim X=\infty$, $\dim Y\ge 1$, and $\dim L^2 = \infty$.
Let $(\g_j)_{j\ge 1}$ be a sequence of independent standard Gaussian random
variables. The
following assertions are equivalent:
\begin{enumerate}

\item[(i)] For all $f\in \Lipo$ and all simple functions $\phi:S\to X$ we have
\[\|u_{f(\phi)}\|_{\ell(L^2, Y)}\leq K \n f\n_{\rm Lip} \|u_{\phi}\|_{\ell(L^2, X)},\]
where $K$ is a constant depending on $X$ and $Y$ only.

\item[(ii)] $X$ has cotype $2$ and $Y$ has type $2$.
\end{enumerate}
If ${\rm (i)}$ holds, then
$C_2^\g(X)\le K$ and $T_2^\g(Y)\le \sqrt{2}K,$
and for
all $f\in\Lipo$ the  mapping $\tilde f$ uniquely extends to an element of
$\text{\rm Lip}_0(\ell(L^2,X),\ell(L^2, Y))$ satisfying $$\n{\tilde f}\n_{\rm
Lip}\le C_2^\g(X)T_2^\g(Y) \n f\n_{\rm Lip}.$$
\end{theorem}
\begin{proof}
(i)$\Rightarrow$(ii): \ Let $x_1, \ldots, x_n\in X$ be
arbitrary. By the $\sigma$-finiteness of $(S,\Sigma,\mu)$ and the assumption
that $\dim L^2=\infty$ there exist disjoint sets $S_1,\dots, S_n\in\Sigma$
satisfying $0<\mu(S_j)<\infty$ for $j=1,\dots,n$ and define $\phi:S\to X$ by
$\phi := \sum_{j=1}^n h_j \otimes x_j$, where $h_j = \mu(S_j)^{-1/2}
\one_{S_j}$ for all $j$. It follows from \eqref{lnorm} that for all $f\in \Lipo$,
\[\begin{aligned}
\Bigl(&\E\Big\|\sum_{j=1}^n \g_j \mu(S_j)^{\frac12} f\big(\mu(S_j)^{-\frac12}
x_j\big)\Big\|^2\Bigr)^{\frac12} = \|u_{f(\phi)}\|_{\ell(L^2, Y)}
\\ & \qquad \leq K
\n f\n_{\rm Lip} \|u_{\phi}\|_{\ell(L^2, X)} = K \n f\n_{\rm Lip}
\Bigl(\E\Big\|\sum_{j=1}^n \g_j x_j\Big\|^2\Bigr)^{\frac12}.
\end{aligned}\]
By an application of Theorem \ref{thm:main} with $a_j =
\mu(S_j)^{-\frac12}$ we obtain (ii).

(ii)$\Rightarrow$(i): Assume that $X$ has cotype $2$ and $Y$ has type $2$ and
fix $f\in\Lipo$. For simple functions $\phi, \psi:S\to X$ we have, by
\eqref{type2} and its cotype $2$ analogue,
$$
\begin{aligned}
\ & \|\tilde{f}(u_\phi) - \tilde{f}(u_{\psi})\|_{\ell(L^2, Y)}
\\ & \qquad   = \| u_{f(\phi)}-u_{f(\psi)}\|_{\ell(L^2, Y)} \leq T_2^\g(Y)\|f(\phi)-f(\psi)\|_{L^2(Y)}
\\ & \qquad \leq T_2^\g(Y) \n f\n_{\rm Lip} \|\phi - \psi\|_{L^2(X)} \leq C_2^\g(X) T_2^\g(Y)\n f\n_{\rm Lip} \|u_\phi -
u_\psi\|_{\ell(L^2, X)}.
\end{aligned}
$$
Since $\ell_{\rm simple}(L^2, X)$ is dense in $\ell(L^2, X)$ it follows that
$\tilde{f}$ has a unique Lipschitz continuous extension from $\ell(L^2, X)$ to
$\ell(L^2, Y)$ with $\n{\tilde{f}}\n_{\rm Lip} \leq  C_2^\g(X)T_2^\g(Y) \n f\n_{\rm Lip}.$ This proves the final assertion, and (i) follows by taking $\psi=0$.
\end{proof}

Theorem \ref{thm:main3} is motivated by the result
from \cite{NW, RS} that a function
$\phi:(0,T)\to X$ is stochastically integrable with respect to a Brownian
motion if and only if the operator $u_\phi$ is well defined and belongs to
$\ell(L^2(0,T),X)$. The question whether $\tilde f$ extends
continuously to $\ell(L^2(0,T),X)$ for all $f\in {\rm Lip}_0(X,X)$ thus
amounts to asking whether $f(\phi)$ is stochastically integrable
whenever $\phi$ has this property.
This question arises naturally in the study of stochastic differential
equations in $X$ driven by multiplicative noise satisfying Lipschitz
conditions; cf. \cite{DZ} for the Hilbert space case.
Theorem \ref{thm:main3} applied to $X=Y$
shows that in general the answer is negative unless $X$
is isomorphic to a Hilbert space.

\medskip
{\em Acknowledgment} -- We thank Tuomas Hyt\"onen and
the anonymous referee for helpful
remarks.

\end{document}